\documentclass[letterpaper, 10 pt, conference]{ieeeconf}

\usepackage{amssymb,amsmath,color}
\usepackage{graphicx}
\usepackage{caption}
\usepackage{enumerate}
\IEEEoverridecommandlockouts                              
\newcommand*{\QEDB}{\hfill\ensuremath{\square}}%

\newtheorem{assumption}{Assumption}
\newtheorem{definition}{Definition}

\newtheorem{lemma}{Lemma}
\newtheorem{theorem}{Theorem}

\newtheorem{proposition}{Proposition}
\newtheorem{remark}{Remark}

\newtheorem{example}{Example}

\begin{document}

\title{\bf Fixed-Time ISS for Interconnected Discrete-Time Systems, Part I:\\ Extinction-Based Lyapunov Functions}

\author{Sumadhu Rubaiyat \thanks{S. Rubaiyat and J. I. Poveda are with the ECE Department, UC, San Diego. Email:\{poveda@ucsd.edu.\}. Research supported in part by NSF grants CMMI 2228791 and CAREER
2305756.}and Jorge I. Poveda.}

\maketitle
\begin{abstract}
We study fixed-time input-to-state stability (FxT-ISS) for discrete-time systems subject to exogenous inputs. We introduce Lyapunov-based sufficient conditions for FxT-ISS based on a class of dissipation functions with a fixed-time extinction property, which encompasses commonly used two-power constructions while allowing substantially more general decay mechanisms. Building on this framework, we establish a small-gain theorem for interconnected FxT-ISS systems in the presence of disturbances. The results are illustrated through several examples, including perturbed optimization algorithms and a discrete-time feedback-optimization architecture.
\end{abstract}
\section{INTRODUCTION}
Traditionally, discrete-time control, optimization, and decision-making algorithms have been analyzed using tools that certify \emph{asymptotic stability} \cite{jiang2001input,Kellet_inclusions1}. This property combines stability and asymptotic attractivity: the former prevents the iterates from undergoing arbitrarily large excursions away from their nominal operating point, while the latter guarantees convergence to that point as time tends to infinity. 
However, increasing computational and performance requirements have motivated the development of algorithms for which convergence is not merely asymptotic, but occurs exactly within a uniformly bounded time. This behavior is captured by the notion of \emph{fixed-time} (FxT) stability. In the continuous-time setting, since its discovery in \cite{andrieu2008homogeneous}, fixed-time stability has received significant attention and has led to new algorithmic architectures for control \cite{lopez2019conditions,polyakov2011nonlinear}, estimation \cite{rios2017time}, learning \cite{poveda2021nonsmooth}, and optimization \cite{garg2022fixed}.
In contrast, fixed-time stability remains considerably less developed in the \emph{discrete-time} setting. A formalization of this property, together with Lyapunov-based certificates for a class of discrete-time systems, was introduced in \cite{lee2023fixed}, with related results subsequently extended to certain stochastic systems \cite{lee2024fixed,tatari2023deterministic}. More recently, the observations in \cite{brogliato2023comments} showed that fixed-time stable discrete-time dynamics can also arise from suitable discretizations of continuous-time controllers. However, despite these developments, several fundamental questions remain open. Existing results have largely focused on particular classes of discrete-time dynamics built around the two-power structures commonly used in continuous-time fixed-time stability \cite{lee2023fixed,tatari2023deterministic}, and more general Lyapunov characterizations remain comparatively unexplored. More importantly, the effect of exogenous inputs and disturbances has received little attention, and existing results have focused on constructive tools for sampled-data controllers in time-varying linear systems with disturbance bounds depending only on a fixed window of recent inputs \cite{mazenc2024fixed}, or ultimate boundedness results \cite{tatari2023deterministic}.  Likewise, systematic conditions guaranteeing preservation of fixed-time stability under interconnection, analogous to fixed-time small-gain results in continuous time \cite{tang2026lyapunov,zimenko2026exponential}, or finite-time results in discrete-time \cite{mo2024lyapunov} are still largely unavailable.

In this paper, we develop a framework for the analysis of fixed-time stability in discrete-time systems subject to exogenous inputs, with particular emphasis on interconnected systems. First, we introduce the notion of \emph{fixed-time input-to-state stability} (FxT-ISS) for discrete-time systems with respect to sets using (generalized) $\mathcal{KL}$ functions, and establish Lyapunov-based sufficient conditions to certify this property. A key technical ingredient is a new class of Lyapunov dissipation functions characterized by a \emph{fixed-time extinction property}: rather than imposing a particular algebraic form on the Lyapunov decrease, this property requires the associated scalar comparison dynamics to vanish after a uniformly bounded number of iterations. This viewpoint provides a unifying framework that recovers the two-power and exponential decay mechanisms studied in \cite{lee2023fixed} and \cite{tatari2023deterministic}, while substantially enlarging the class of admissible dissipation functions and providing a more flexible mechanism for certifying fixed-time convergence in discrete time, enabling the uniform incorporation of exogenous disturbances in the analysis via ISS tools. Building on this characterization, and inspired by recent results in continuous-time \cite{tang2026lyapunov,tang2026smallgain}, we then develop a \emph{fixed-time small-gain-type theorem} for interconnected discrete-time systems whose individual components satisfy suitable FxT-ISS properties. The result provides conditions under which the fixed-time behavior of the individual subsystems is preserved under feedback interconnection, while explicitly accounting for exogenous disturbances. Together, these results provide one particular extension of the Lyapunov and small-gain machinery commonly available for asymptotically stable discrete-time systems to the fixed-time regime. For alternative Lyapunov-based characterizations suitable for FxT-ISS systems that do not admit the Lyapunov dissipation functions studied in this paper, we refer the reader to the companion paper \cite{ACC26B}. 


\section{PRELIMINARIES}
\label{sec_notation}
\subsection{Notation}
Given a nonempty closed set $\mathcal{A} \subset \mathbb{R}^n$ and a vector $z \in \mathbb{R}^n$, we define $|z|_\mathcal{A} := \min_{y \in \mathcal{A}}|z-y|$, and we use $|\cdot|$ to denote the standard Euclidean norm.  Throughout the paper, for ease of notation we write $(u,v)$ for $(u^{T},v^{T})^{T}$.  We let $\mathcal{PD}$ be the class of continuous functions $\rho:\mathbb{R}_{\geq0}\to\mathbb{R}_{\geq0}$ that are positive definite. A continuous function $\alpha:\mathbb{R}_{\geq0}\to\mathbb{R}$ is said to be of class $\mathcal{K}$ if it is zero at zero and strictly increasing, and it is of class $\mathcal{K}_{\infty}$ if, additionally, it is unbounded. A function $\beta:\mathbb{R}_{\geq 0}\times\mathbb{Z}_{\geq 0}
\to\mathbb{R}_{\geq 0}$ belongs to $\mathcal{KL}_{\mathrm{FxT}}$ if there exists
$M\in\mathbb{Z}_{\geq 1}$ such that: (a) for every $k\in\mathbb{Z}_{\geq 0}$, $\beta(\cdot,k)$ is continuous, nondecreasing, and $\beta(0,k)=0$; (b) for every $k<M$, $\beta(\cdot,k)$ belongs to class $\mathcal{K}$; (c) for every $r\geq0$, $\beta(r,\cdot)$ is  nonincreasing; (d) $\beta(r,k)=0$ for every $r\geq0$ and every $k\geq M$. For any $r\in\mathbb{R} $, we use $[r]_+:=\max\{r,0\}$, and for any map $h:\mathbb R_{\geq0}\to\mathbb R_{\geq0}$,
we use $h^k$ to denote its $k$-fold composition, defined
recursively by $h^0:=\operatorname{id}$, $h^{k+1}:=h\circ h^k$.

\begin{table}[t]
\centering
\label{tab:extinction_functions}
\scriptsize
\setlength{\tabcolsep}{2pt}
\renewcommand{\arraystretch}{1.15}
\resizebox{\columnwidth}{!}{%
\begin{tabular}{|ll@{\hspace{6pt}}ll|}
\hline
$\psi(s)$ & Conditions & $\psi(s)$ & Conditions\\
\hline
$s$ & ---
& $s+cs^p-r(1-e^{-s})$ & $r\le1$\\

$(as^p+bs^q)^r$
& $p,q,r>0,\;pr<1<qr$
& $s-c\frac{t}{1+t}\sin^2s$ & $c\le\varepsilon$\\

$a\max\{s^p,s^q\}$ & $q>1$
& $s-c\min\{t,1\}$ & $c\le\varepsilon$\\

$as^p(1+s)^q$ & $q>1-p$
& $s-c(1-e^{-t})$ & $c\le\varepsilon$\\

$as^p e^{bs^q}$ & $q>0$
& $s-c\tanh t$ & $c\le\varepsilon$\\

$\min\{\kappa^{-1}e^s,s\}$ & ---
& $s-c\frac{t^p}{1+t^p}$ & $p>0,\;c\le\varepsilon$\\

$a(e^{s^p}-1)$ & ---
& $s-c\frac{\log(1+t)}{1+\log(1+t)}$
& $c\le\varepsilon$\\

$a\sinh(s^p)$ & ---
& $s-c\,t e^{1-t}$ & $c\le\varepsilon$\\

$a\log(1+s^p)+bs^q$ & $q>1$
& $s-c[1-\lvert t-1\rvert]_+$ & $c\le\varepsilon$\\

$\frac{as^p}{1+s^p}+\frac{s^2}{1+s}$ & ---
& $s-\min\{t,L\}$ & ---\\

$\max\{as^p,s-H\}$ & ---
& & \\
\hline
\end{tabular}%
}
\caption{Examples of functions $\psi\in\mathcal{E}_{\text{FxT}}^M$.
All parameters are positive, and $0<p<1$ unless otherwise
specified. Here $t=[s-\varepsilon]_+$. A dash indicates no further restrictions.}
\vspace{-0.5cm}
\end{table}
\subsection{The Fixed-Time Extinction Property}
The following definition will be used throughout the paper.

\vspace{0.1cm}
\begin{definition}
A function $\psi\in\mathcal{PD}$ is said to have the \emph{fixed-time extinction property} if there exists $M\in\mathbb{Z}_{\geq1}$ such that the function
\begin{equation}
    h_{\psi}(s):=
    \max_{0\leq r\leq s}g(r),~~s\geq0,~~\text{with}~~g(r):=[r-\psi(r)]_{+},
    \label{eq:general_h}
\end{equation}
satisfies $h_{\psi}^M(s)=0$ for all $s\geq0$. In this case, we say that $\psi\in\mathcal{E}^M_{\text{FxT}}$. \QEDB
\end{definition}

\vspace{0.1cm}
The following lemma provides an equivalent characterization of the set $\mathcal{E}_{\text{FxT}}^M$.
 
\vspace{0.1cm}
\begin{lemma} \label{lem:general_fxt_extinction}
 Let $\psi\in\mathcal{PD}$. Then there exists
$M\in\mathbb Z_{\geq1}$ such that
$\psi\in\mathcal E_{\mathrm{FxT}}^M$ if and only if there exist $\varepsilon>0$ and $H<\infty$ such that: (a) $\psi(s)\geq s$ for all $0\leq s\leq\varepsilon$; (b) $\psi(s)\geq s-H$ for all $s\geq0$. Moreover, if $\psi\in\mathcal E_{\mathrm{FxT}}^M$
for some $M\in\mathbb Z_{\geq1}$, then
$h_\psi(s)<s$ for every $s>0$.
\hfill $\square$
\end{lemma}

\vspace{0.1cm}
Table I shows different examples of functions $\psi\in\mathcal E_{\mathrm{FxT}}^M$. Discrete-time analyses of FxT stability that involve the functions $\psi(s)=as^\alpha+bs^\beta$ and $\psi(s)=\kappa_f^{-1}e^{s^p}s^{1-p}$ appeared in \cite[Thm. 3.1 $\&$ 4.1]{lee2023fixed}, and $\psi(s)=\alpha\max\{s^{r_1},s^{r_2}\}$  in \cite[Thm. 1]{tatari2023deterministic}. However, as shown in the next section, FxT stability admits a much broader class of dissipation functions $\psi$.

%
%
\vspace{0.1cm}

\section{EXTINCTION-BASED LYAPUNOV CONDITIONS FOR FXT-ISS}
\label{mainresults}
In this paper, we consider time-invariant discrete-time dynamical systems with inputs, given by
\begin{equation}\label{eq:dynamics}
x^+=g(x,u),~~~~x_0:=x(0)\in\mathbb{R}^n,
\end{equation}
where $g:\mathbb{R}^n\times\mathbb{R}^m\to\mathbb{R}^n$ is the state-update map, and $u:\mathbb{Z}_{\geq0}\to\mathbb{R}^m$ is an input signal. For these systems, we are interested in studying stability properties relative to general closed sets $\mathcal{A}\subset\mathbb{R}^n$. For this, we make the following regularity assumption.
\begin{assumption}\label{assumption1}
The function $g(\cdot,\cdot)$ is continuous, and the set $\mathcal{A}$ is nonempty and closed. \hfill $\QEDB$
\end{assumption}

The following definition formalizes the property of \emph{fixed-time input-to-state stability}.

\vspace{0.1cm}
\begin{definition}\label{def:fxISS}
Let $\mathcal{A}\subset\mathbb{R}^n$ be a nonempty closed set. System \eqref{eq:dynamics} is said to render the set $\mathcal{A}$ fixed-time input-to-state stable (ISS) if there exists a class $\mathcal{K}\mathcal{L}_{\text{FxT}}$ function $\beta$ and a class $\mathcal{K}$ function $\gamma(\cdot)$, such that for all $x_0\in\mathbb{R}^n$ and all inputs $u\in\ell^{\infty}$, all solutions of \eqref{eq:dynamics} satisfy the following inequality
\begin{equation}\label{eq:general_fxtiss_estimate}
    |x(j)|_{\mathcal A}\leq\max\left\{\beta(|x(0)|_{\mathcal A},j),\gamma\left(|u|_{[0,j-1],\infty}\right)\right\},
\end{equation}
for all $j\geq0$, where $|u|_{[0,j-1],\infty}=0$ for $j=0$.  \hfill $\square$
\end{definition}

Definition \eqref{def:fxISS} is a natural extension to the discrete-time setting of similar FxT-ISS notions in continuous time \cite{lopez2019conditions}. 
\subsection{Sufficient Lyapunov-based Conditions}
The following definition introduces a particular type of ISS-Lyapunov functions for fixed-time stability in discrete time.
\vspace{0.1cm}
\begin{definition} \label{def:general_fxtiss_lyap}
A continuous function $V:\mathbb R^n\to\mathbb R_{\geq0}$ is called an \emph{Extinction-based FxT-ISS
Lyapunov function} with respect to the nonempty closed set $\mathcal A\subset\mathbb{R}^n$ if there exist $\alpha_1,\alpha_2\in\mathcal K_\infty$, $\rho\in\mathcal K_\infty$, and $\psi\in\mathcal{E}^M_{\text{FxT}}$, such that:
\begin{subequations}\label{eq:suffLyap1}
\begin{align}   
    \alpha_1(|x|_{\mathcal A})&\leq V(x)\leq \alpha_2(|x|_{\mathcal A}),
    \qquad \forall x\in\mathbb R^n, \label{eq:general_lyap_bounds}\\
    V(x^+)-V(x)&\leq-\min\{\psi(V(x)),V(x)\}+\rho(|u|)\label{eq:additive_dissipation}
\end{align}
\end{subequations}
for every $(x,u)\in\mathbb R^n\times\mathbb R^m$. \hfill $\square$
\end{definition}

\vspace{0.1cm}
Under Assumption \ref{assumption1} and compactness of $\mathcal{A}$, the following implication-form condition provides a sufficient condition to obtain an Extinction-Based FxT-ISS Lyapunov function for system \eqref{eq:dynamics}:

\vspace{0.1cm}
\begin{proposition}\label{prop:implication_to_dissipation}
Consider the discrete-time system \eqref{eq:dynamics}. Suppose that 
Assumption \ref{assumption1} holds and that $\mathcal{A}$ is compact.  Let $V:\mathbb{R}^n\to\mathbb{R}_{\geq0}$ be continuous and satisfy conditions \eqref{eq:general_lyap_bounds} with $\alpha_1,\alpha_2\in\mathcal{K}_{\infty}$, and suppose there exist $\phi\in\mathcal{K}_{\infty}$ and  $\psi\in\mathcal{E}_{FxT}^M$ such that for all $(x,u)\in\mathbb{R}^n\times\mathbb{R}^m$: 
\begin{align}  \label{eq:prop_implication}
V(x)\geq&~\phi(|u|)\implies \notag\\
&V(x^+)-V(x)\leq-\min\{\psi(V(x)),V(x)\}
\end{align}
Then, $V$ is an Extinction-Based FxT-ISS Lyapunov function for system \eqref{eq:dynamics}.  \hfill $\QEDB$
\end{proposition}

\vspace{0.1cm}
Unlike the standard (asymptotic) ISS case \cite{jiang2001input,cai2008input}, under Assumption \ref{assumption1} and compactness of $\mathcal{A}$ the sufficient condition of Proposition \ref{prop:implication_to_dissipation} is not equivalent to Definition \ref{def:general_fxtiss_lyap}. We illustrate this fact with the following stylized example.

\vspace{0.1cm}
\begin{example}[FxT-ISS of the Perturbed Newton Method]
\label{ex:disturbed_newton}
Consider the quadratic function $f(x)=\frac{1}{2}x^\top Qx+b^\top x+c$, where $Q=Q^\top\succ 0$. Its unique minimizer is $x^\star=-Q^{-1}b$. Consider the following perturbed Newton optimization algorithm:
\begin{equation}
    x^+ =x-\bigl[\nabla^2 f(x)\bigr]^{-1}
        \nabla f(x)+u,
    \label{eq:disturbed_newton}
\end{equation}
where $u\in\mathbb R^n$ represents an additive disturbance. For this system, the function $V(x):=\|x-x^\star\|$ is a FxT-ISS Lyapunov function with respect to $\mathcal A=\{x^\star\}$. Indeed, since $\nabla f(x)=Qx+b$ and $\nabla^2 f(x)=Q,$ the iteration \eqref{eq:disturbed_newton}
reduces to $x^+=x-Q^{-1}(Qx+b)+u=-Q^{-1}b+u=x^\star+u$. Consequently, the error $e:=x-x^\star$ satisfies $e^+=u$. By construction, $|x|_{\mathcal A}=V(x)=\|x-x^\star\|$. Thus, condition \eqref{eq:general_lyap_bounds} holds with $\alpha_1(s)=\alpha_2(s)=s$. Moreover, $ V(x^+)=|u|$, and we get $V(x^+)-V(x)= -V(x)+|u|= -\min\{\psi(V(x)),V(x)\}+\rho(|u|)$, where $\psi(s)=s$ and $ \rho(s)=s$.  Thus, $V$ is a Extinction-Based FxT-ISS Lyapunov function for \eqref{eq:disturbed_newton}.  However, the following lemma shows that \eqref{eq:disturbed_newton} does not admit a function $V$ that satisfies the conditions of Proposition \ref{prop:implication_to_dissipation}.  
\begin{lemma}\label{prop:newton_implication_impossibility}
Let $Q=Q^\top\succ0$. Then, for the perturbed Newton iteration \eqref{eq:disturbed_newton} there is no function
$V:\mathbb R^n\to\mathbb R_{\geq0}$
(not necessarily continuous), functions $\alpha_1,\alpha_2,\phi\in\mathcal K_\infty$,
and a function $\psi\in\mathcal E_{\mathrm{FxT}}^M$, for any finite $M\geq1$, such that inequalities \eqref{eq:general_lyap_bounds} and \eqref{eq:prop_implication} hold with $\mathcal{A}=\{x^*\}$ for every $(x,u)\in\mathbb R^n\times\mathbb R^n$. \QEDB 
\end{lemma}
\end{example}

\vspace{0.1cm}
To study systems with discontinuous dynamics $g$, we can relax Assumption \ref{assumption1} in Proposition \ref{prop:implication_to_dissipation},  at the expense of an additional growth condition on $V$:
\begin{proposition}\label{prop:implication_fxt_iss}
Let $\mathcal{A}$ be a nonempty closed set, and consider the discrete-time system \eqref{eq:dynamics}. Let $V:\mathbb R^n\to\mathbb R_{\geq 0}$
be continuous and satisfy \eqref{eq:general_lyap_bounds}
with $\alpha_1,\alpha_2\in\mathcal K_\infty$.
Suppose there exist
$\phi,\kappa\in\mathcal K_\infty$ and
$\psi\in\mathcal E_{\mathrm{FxT}}^M$
such that, for all
$(x,u)\in\mathbb R^n\times\mathbb R^m$,
\begin{subequations}
\label{eq:implication_conditions}
\begin{align}
V(x)\geq&~\phi(|u|)
\implies\notag\\
&V(x^+)-V(x)
\leq-\min\{\psi(V(x)),V(x)\},\label{eq:implication_condition}\\
V(x)\leq&~\phi(|u|)
\implies
V(x^+)\leq\kappa(|u|).
\label{eq:subthreshold_condition}
\end{align}
\end{subequations}
Then $V$ is an Extinction-Based FxT-ISS Lyapunov function
for system \eqref{eq:dynamics}. \hfill $\square$
\end{proposition}

\vspace{0.1cm}
The following theorem, which is the first main result of this paper, shows that the existence of a FxT-ISS Lyapunov function is a sufficient condition for fixed-time ISS.
\begin{theorem}\label{thm:general_fxtiss}
Let $\mathcal{A}\subset\mathbb{R}^n$ be a nonempty closed set, and suppose that $V$ is a FxT-ISS Lyapunov function for system  \eqref{eq:dynamics}. Then, $\mathcal{A}$ is FxT-ISS for system \eqref{eq:dynamics}. \QEDB
\end{theorem}

\vspace{0.1cm}
To the best of our knowledge, the result of Theorem \ref{thm:general_fxtiss} is the first Lyapunov-based sufficient condition of the FxT-ISS property in discrete-time systems. For the case when $u=0$, Theorem \ref{thm:general_fxtiss} generalizes existing FxT stability results for discrete-time systems, recovering \cite{lee2023fixed,tatari2023deterministic} for particular choices of $\psi$, see Table 1. However, by Lemma 1, the class of functions $\psi$ that can be used in Theorem \ref{thm:general_fxtiss} is much broader. Moreover, a FxT-ISS Lyapunov function can be obtained via the inequalities in Propositions \ref{prop:implication_to_dissipation} or \ref{prop:implication_fxt_iss}.

\vspace{0.1cm}
\begin{remark}
By Theorem \ref{thm:general_fxtiss}, the perturbed Newton-method \eqref{eq:disturbed_newton} is FxT-ISS.  While the fixed-time (one-step) convergence property of this method is trivial to prove, the uniform $\mathcal{KL}_{\text{FxT}}$ bound \eqref{eq:general_fxtiss_estimate} provides an additional stability characterization and quantifies the effect of $u$ on the dynamics, uniformly across initial conditions and inputs. Moreover, as shown in the next sections, it opens the door for the study of FxT-ISS in other algorithms, as well as interconnected systems via small-gain arguments.
\QEDB
\end{remark} 
\subsection{Examples and Applications of FxT-ISS}
We present several examples to illustrate the results.
\subsubsection{FxT-ISS of projected Gradient Descent}
\label{ex:projected_gradient_fxt_iss}
Let $C\subset\mathbb R^n$ be a nonempty, closed,
and convex set, and consider the constrained
optimization problem
\begin{equation}
    \min_{x\in C} f(x),
    \label{eq:projected_gradient_problem}
\end{equation}
where $f(x):=\frac{1}{2}\|x-a\|^2$, and $a\in\mathbb R^n$. The unique minimizer is given by $x^\star=\Pi_C(a)$, where $\Pi_C$ denotes the Euclidean projection onto $C$. To find $x^\star$, consider the following perturbed projected-gradient algorithm
\begin{equation}
    x^+
    =\Pi_C\bigl(x-\nabla f(x)+u\bigr),
    \qquad u\in\mathbb R^n.
    \label{eq:disturbed_projected_gradient}
\end{equation}
We show that this system is FxT-ISS with respect
to $\mathcal A:=\{x^\star\}$. Indeed, since $\nabla f(x)=x-a$, the dynamics reduce to $x^+=\Pi_C(a+u)$. To study this system, consider the Lyapunov function $V(x):=\|x-x^\star\|$. Clearly, the comparison bounds \eqref{eq:general_lyap_bounds} hold with $\alpha_1(s)=\alpha_2(s)=s$. Moreover, by the nonexpansiveness of the Euclidean projection onto a closed convex set, we have $V(x^+)=\|\Pi_C(a+u)-\Pi_C(a)\|\leq|u|$. Consequently,
\begin{equation}
\begin{aligned}
    V(x^+)-V(x)
    &\leq -V(x)+|u|\\
    &=-\min\{\psi(V(x)),V(x)\}
      +\rho(|u|),
\end{aligned}
    \label{eq:projected_gradient_dissipation}
\end{equation}
where $\psi(s):=s$ and $ \rho(s):=s$. The function $\rho$ belongs to
$\mathcal K_\infty$. Therefore, $V$ is a Extinction-Based FxT-ISS Lyapunov function and by Theorem \ref{thm:general_fxtiss} the point $x^{\star}$ is FxT-ISS for \eqref{eq:disturbed_projected_gradient}. \QEDB 

\vspace{0.1cm}
\begin{figure}[t]
    \centering
    \includegraphics[width=\columnwidth]{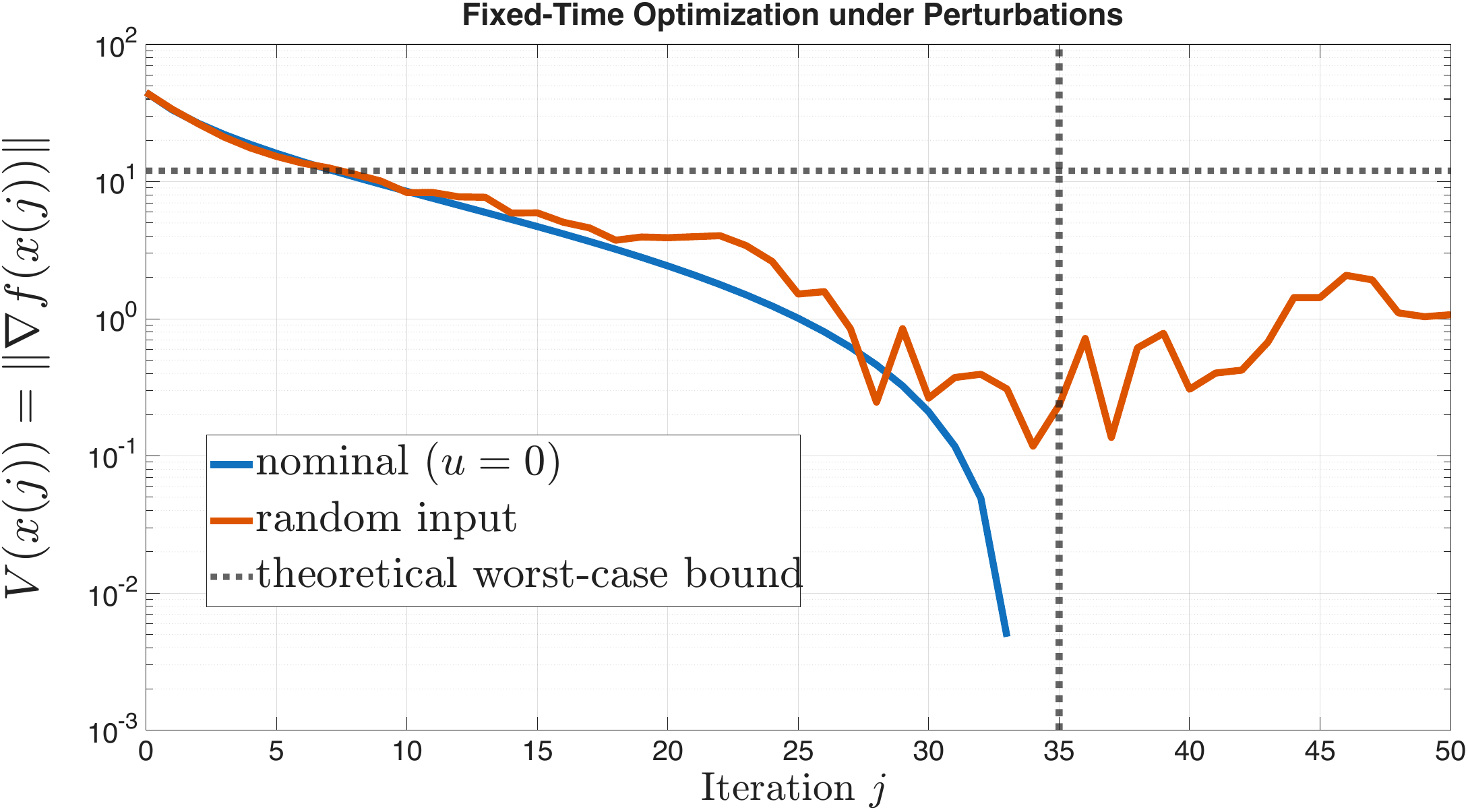}
    \caption{Evolution of $V$ under system \eqref{eq:inverse_gradient_algorithm}. The vertical line shows the theoretical upper bound on the convergence time.}
    \label{fig:feedback_optimization}
    \vspace{-0.5cm}
\end{figure}
\subsubsection{Fixed-time Optimization of Strongly Convex Functions under Inverse Gradient Oracle}
Consider the optimization problem \eqref{eq:projected_gradient_problem} with $C=\mathbb{R}^n$, where $f:\mathbb R^n\to\mathbb R$ is continuously differentiable, $\mu$-strongly convex, has a globally
$L$-Lipschitz gradient, and $x^\star$ is its unique minimizer. Let $G(\cdot):=\nabla f(\cdot)$. Strong convexity implies that $G$ is bijective and that its inverse $G^{-1}:\mathbb R^n
\to\mathbb R^n$ is globally Lipschitz. Let $a,b>0$,  $0<\alpha<1<\beta$, and $\psi$ be given as in Lemma 2-(b). We consider the (perturbed) optimization algorithm
\begin{equation}
    x^+=g(x,u):=G^{-1}\bigl(\mathcal R(G(x))\bigr)+u,
    \label{eq:inverse_gradient_algorithm}
\end{equation}
where $u\in\mathbb R^n$ is an external disturbance and
\begin{equation}
    \mathcal R(z):=
    \begin{cases}
    \displaystyle
    \frac{[|z|-\psi(|z|)]_+}{|z|}\,z,
        & z\neq 0,\\[2mm]
    0, & z=0.
    \end{cases}
    \label{eq:gradient_radial_map}
\end{equation}
We show that \eqref{eq:inverse_gradient_algorithm} admits a FxT-ISS Lyapunov
function with respect to $\mathcal A=\{x^\star\}$. Indeed, let $V(x):=|\nabla f(x)|=|G(x)|$ and note that by strong convexity and Lipschitz continuity of the gradient we have $\mu|x-x^\star|\le V(x)\le L|x-x^\star|$. Thus, \eqref{eq:general_lyap_bounds} holds with
$\alpha_1(s)=\mu s$ and $\alpha_2(s)=Ls$. To establish \eqref{eq:suffLyap1}, let
$z=G(x)$ and $\bar x=G^{-1}(\mathcal R(z))$.
By Lipschitz continuity of $G$,
\begin{align}
    V(x^+)&=|G(\bar x+u)|\le |G(\bar x)|+L|u|=|\mathcal R(z)|+L|u| \nonumber\\
    &=[V(x)-\psi(V(x))]_++L|u|.
    \label{eq:inverse_gradient_bound}
\end{align}
Subtracting $V(x)$ from both sides establishes \eqref{eq:additive_dissipation}, Therefore, $V$ is an Extinction-Based FxT-ISS Lyapunov function for \eqref{eq:inverse_gradient_algorithm},
with linear gain $\rho(s)=Ls$. By Theorem \ref{thm:general_fxtiss}, $x^{\star}$ is FxT-ISS. 

\vspace{0.1cm}
\subsubsection{Planar System with Radial State-Feedback}
Let $x \in \mathbb{R}^2$, and consider the open-loop system
\begin{equation}
x^+= (1+v)R(\theta)x,
\end{equation}
where $\theta \in \mathbb{R}$, $R(\theta)$ is a $2 \times 2$ orthogonal rotation matrix (i.e., $\Vert{}R(\theta)\Vert{} = 1$)) and $v$ is the control signal. When $v=0$, this system simply rotates the state $x$ without changing its magnitude. To stabilize the origin in fixed time, consider first the unsaturated feedback law: $\hat{v}(x,u)=- c_2 |x|^{a-1}- d_2 |x|^{\beta-1}+\gamma u$, where $\gamma>0$, $a \in (0, 1)$, $\beta \in (1, \infty)$, $c_2 >0$, $d_2 >0$ and $\gamma u$ is a disturbance.  Let the feedback law be 
\begin{equation}
v(x,u)=
\begin{cases}
\Big[1+\hat{v}(x,u)\Big]_+-1,
& x\neq 0,\\
-1,
& x=0,
\end{cases}
\end{equation}
which imposes a lower saturation on the feedback law, leading to the following closed-loop system:
\begin{equation*}\label{cldexample}
x^+=\left\{\begin{array}{ll} 
\left[1 - c_2 |x|^{a-1} - d_2 |x|^{\beta-1}+\gamma u\right]_+\cdot R(\theta) x &x\neq0\\
0 & x=0,
\end{array}\right.
\end{equation*}
where  $\theta \in \mathbb{R}$, $R(\theta)$ is a $2 \times 2$ orthogonal rotation matrix (i.e., $\Vert{}R(\theta)\Vert{} = 1$)). Note that the right-hand side of this system is continuous in $(x,u)$. 

To study the stability properties of the closed-loop dynamics, we use Proposition \ref{prop:implication_to_dissipation} with $V(x)=|x|$. Since $R(\theta)$ is orthogonal, for every $x\neq 0$,
\begin{align}
V(x^+)&=
\Big[V(x)-c_2V(x)^a-d_2V(x)^\beta+\gamma uV(x)\Big]_+\notag\\
&\leq\Big[V(x)-c_2V(x)^a-d_2V(x)^\beta+\gamma |u|V(x)\Big]_+ .\label{eq:lyap_bound}
\end{align}
%
The same inequality holds at $x=0$. Let $c\in(0,\min\{c_2,1\})$, $d\in(0,\min\{d_2,1\})$, and define the class-$\mathcal{K}$ function $\phi(s):=\left(\frac{\gamma s}{d_2-d}\right)^{\frac{1}{\beta-1}}$. Whenever $|x|\geq\phi(|u|)$, it follows that $\gamma|u|V(x)\leq(d_2-d)V(x)^\beta$. Substituting this inequality into \eqref{eq:lyap_bound} yields
\begin{align}
V(x^+)&\leq\Big[V(x)-c_2V(x)^a-dV(x)^\beta\Big]_+\notag\\
&\leq\Big[V(x)-cV(x)^a-dV(x)^\beta
\Big]_+ .
\end{align}
Using the identity $[\tilde{a}-\tilde{b}]_+-\tilde{a}=-\min\{\tilde{a},\tilde{b}\}$, for all $\tilde{a},\tilde{b}\geq 0$ we conclude that $V(x^+)-V(x)\leq-\min\{cV(x)^a+dV(x)^\beta,V(x)\}$. Therefore, $V$ satisfies \eqref{eq:prop_implication}, and by Proposition \ref{prop:implication_to_dissipation}, it is an Extinction-Based FxT-ISS Lyapunov function. By Theorem~\ref{thm:general_fxtiss} the set $\mathcal{A}=\{0\}$ is FxT-ISS. \hfill $\square$
\section{AN EXTINCTION-BASED SMALL-GAIN THEOREM FOR FXT-ISS}
\label{sec:smallgain}
In this section, we leverage the FxT-ISS characterization introduced in the previous section to establish a Fixed-Time small-gain theorem for interconnected discrete-time systems, of the form 
\begin{figure}[t]
    \centering
    \includegraphics[width=0.9\columnwidth]{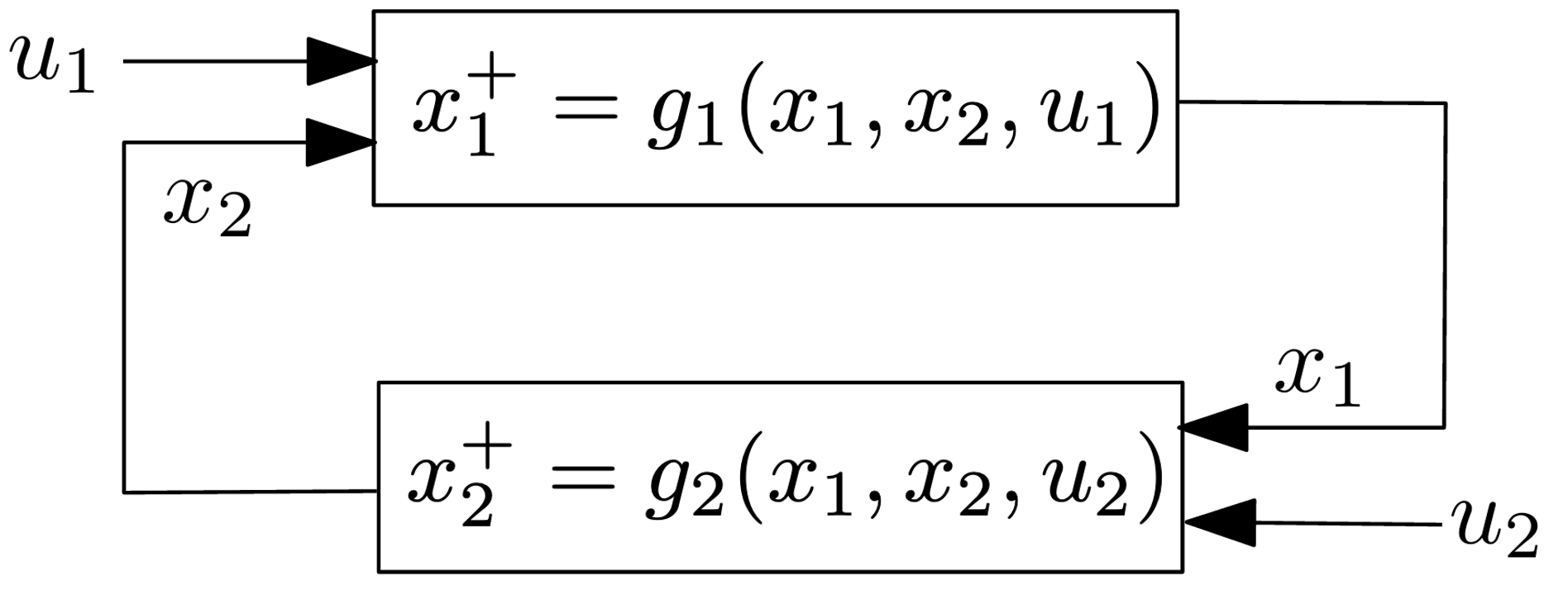}
    \caption{Interconnection between two sub-systems under perturbations.}
    \label{fig:smallgain}
    \vspace{-0.5cm}
\end{figure}
\begin{equation}\label{eq:sg_interconnection}
\begin{aligned}
x_1^+ &= g_1(x_1,x_2,u_1),\\
x_2^+ &= g_2(x_2,x_1,u_2),
\end{aligned}
\end{equation}
where $x_i\in\mathbb R^{n_i}$ and $u_i\in\mathbb R^{m_i}$, and suppose that $g_i(0,0,0)=0$ for $i\in\{1,2\}$.  Figure \ref{fig:smallgain} shows a block diagram representation of this system. 
\subsection{Main Assumptions and Result}
We impose the following assumptions on each of the subsystems. Below, we use $x:=(x_1,x_2)$ and $u:=(u_1,u_2)$.
\begin{assumption}\label{ass:sg_subsystems}
For every $i\in\{1,2\}$ and $j\neq i$, there exist a continuous
function $V_i:\mathbb R^{n_i}\to\mathbb R_{\geq0}$, functions $\underline{\alpha}_i,\overline{\alpha}_i,
\sigma_i,\chi_i\in\mathcal K_\infty$, a function $\psi_i\in\mathcal E_{\mathrm{FxT}}^{M_i}$, and a continuous, nonnegative, nondecreasing function $\gamma_{ij}:\mathbb R_{\geq0}\to\mathbb R_{\geq0}$, with $\gamma_{ij}(0)=0$, such that
\begin{equation}\label{eq:sg_subsystem_bounds}
\underline{\alpha}_i(|x_i|)
\leq V_i(x_i)
\leq\overline{\alpha}_i(|x_i|),~~~\forall~x_i\in\mathbb{R}^{n_i},
\end{equation}
and
\begin{align}
&V_i\!\left(x_i^+\right)-V_i(x_i)
\notag\\
&\quad\leq
\max\Big\{
-\min\{\psi_i(V_i(x_i)),V_i(x_i)\},
\gamma_{ij}(V_j(x_j))\notag\\
&~~~~~~~~~~~~~~~-V_i(x_i)\Big\}+\chi_i(|u_i|),
\label{eq:sg_subsystem_gain}
\end{align}
for all $(x_i,x_j,u_i)$.
\hfill $\square$
\end{assumption}
\begin{remark}
Note that adding $V_i(x_i)$ to both sides of \eqref{eq:sg_subsystem_gain} and using
$s-\min\{\psi_i(s),s\}=[s-\psi_i(s)]_+$ gives the equivalent inequality
\begin{align}
V_i\!\left(x_i^+\right)
\leq
\max\Big\{
&[V_i(x_i)-\psi_i(V_i(x_i))]_+,
\gamma_{ij}(V_j(x_j))
\Big\}\notag\\
&+\chi_i(|u_i|).
\label{eq:sg_subsystem_max_form}
\end{align}
In this sense, Assumption \ref{ass:sg_subsystems} is similar in nature to \cite[Assumption 1]{tang2026smallgain} in the context of FxT-ISS for interconnected continuous-time systems. \hfill $\square$
\end{remark}

To establish a FxT small-gain theorem for system \eqref{eq:sg_interconnection}, we define the following functions:
\begin{subequations}\label{eq:sg_scaled_maps}
\begin{align}
a_i(s)
&:=
\sigma_i^{-1}\circ h_{\psi_i}\circ\sigma_i(s),
\label{eq:sg_scaled_decay}\\
b_{ij}(s)
&:=
\sigma_i^{-1}\circ\gamma_{ij}\circ\sigma_j(s).
\label{eq:sg_scaled_gain}
\end{align}
\end{subequations}
where $h_{\psi}$ was defined in \eqref{eq:general_h}. We impose the following regularity conditions on these functions:
\begin{assumption}\label{ass:sg_fxt_condition}
For every $i\in\{1,2\}$ and $j\neq i$, the scaled internal
gain $b_{ij}$ is globally bounded and vanishes on a neighborhood of the origin. \QEDB
\end{assumption}

By leveraging the previous constructions, we can establish the following theorem, which is the second main result of this paper:
\begin{theorem}\label{thm:gain_fxt_small_gain}
Suppose that Assumptions~\ref{ass:sg_subsystems}
and~\ref{ass:sg_fxt_condition} hold and that, for every
$i\in\{1,2\}$ and $j\neq i$,
\begin{equation}\label{eq:sg_internal_gain_condition}
b_{ij}(s)<s,
\qquad
\forall s>0.
\end{equation}
Then the interconnected system \eqref{eq:sg_interconnection} renders the set $\mathcal{A}=\{0\}$ FxT-ISS with respect to the input $u=(u_1,u_2)$.\hfill $\QEDB$
\end{theorem}

\vspace{0.1cm}
\begin{remark}
Condition \eqref{eq:sg_internal_gain_condition} implies
\[
b_{12}\circ b_{21}(s)<s,
\qquad
b_{21}\circ b_{12}(s)<s,
\qquad s>0,
\]
which is the usual cyclic small-gain condition. The additional boundedness and dead-zone requirements are specific to the discrete-time fixed-time setting, and they ensure exact extinction after a uniformly bounded number of iterations rather than merely asymptotic convergence. \hfill $\square$
\end{remark}
\subsection{Application: Fixed-time Feedback Optimization in Discrete-Time}
\label{ex:fxt_feedback_optimization}
Consider the optimization problem \eqref{eq:projected_gradient_problem} with $C=\mathbb{R}^n$, where $f:\mathbb R^n\to\mathbb R$ is continuously differentiable, $\mu$-strongly convex, has a
globally $L$-Lipschitz gradient, and its unique minimizer is $y^\star$. The argument of $f$ is the state of the discrete-time plant $y^+=v+d_p$, where $y,v,d_p\in\mathbb R^n$ denote the plant state, control input, and external disturbance, respectively. Let $z\in\mathbb R^n$ denote the reference generated by a discrete-time optimization algorithm, and let $d_o\in\mathbb R^n$ represent disturbances
affecting its update. The goal of the optimization algorithm is to regulate $y$ towards $y^{\star}$ while keeping the plant stable. %
\begin{figure}[t]
    \centering
    \includegraphics[width=\columnwidth]{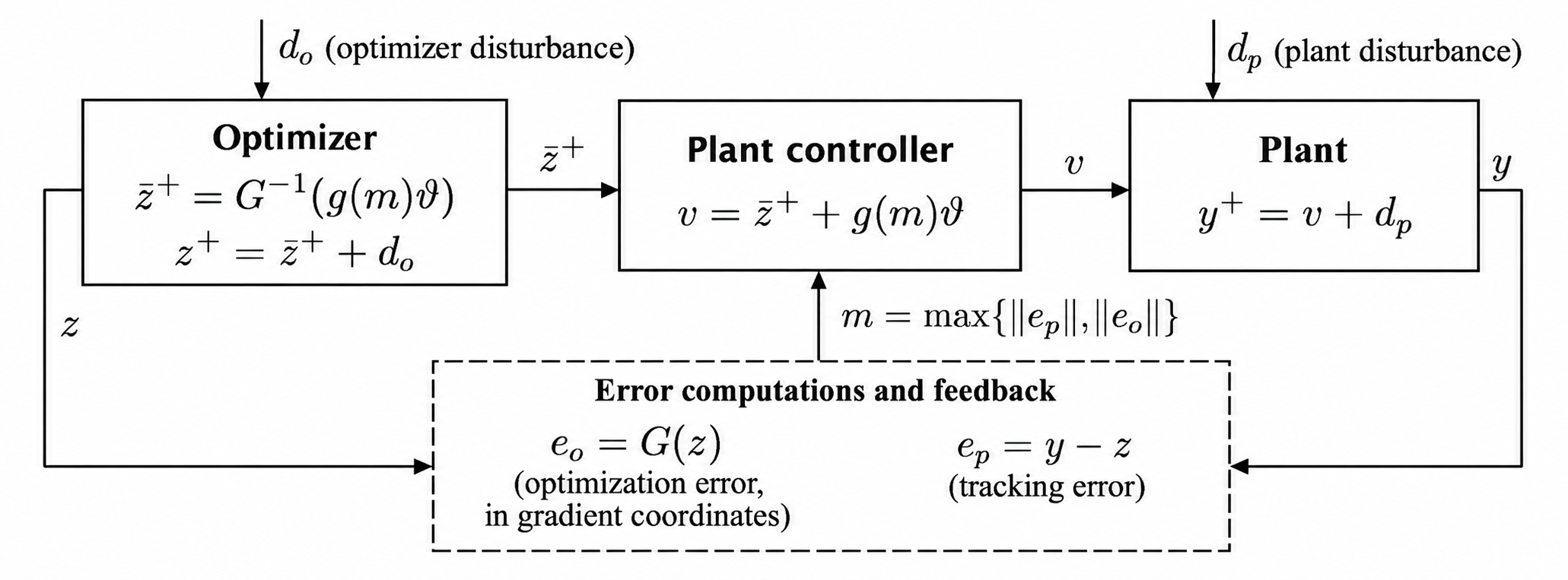}
    \caption{Block diagram of the feedback optimization
    interconnection with plant and optimizer disturbances
    $d_p$ and $d_o$, respectively.}
    \label{fig:feedback_optimization}
    \vspace{-0.5cm}
\end{figure}
The plant tracking error and the optimization error (in gradient coordinates) are given by $e_p:=y-z$ and $e_o:=G(z)$, respectively. Since $G$ is bijective, $e_o=0$ if and only if
$z=y^\star$. Let $g(s):=3\frac{[s-1]_+}{1+[s-1]_+}$, for all $s\ge0$, and note that \(g\) is continuous, nondecreasing, globally bounded by \(3\), and vanishes on \([0,1]\). Moreover, \(g(s)<s\) for all \(s>0\), and \(g^5(s)=0\) for all \(s\ge0\). Fix an arbitrary unit vector
$\vartheta\in\mathbb R^n$ and define the
nominal next optimizer reference
\begin{equation}
    \bar z^+:=G^{-1}\left(g\!\left(\max\{|e_p|,|e_o|\}\right)\vartheta\right).
    \label{eq:fo_nominal_reference}
\end{equation}
Based on this,  we consider the perturbed optimization algorithm $z^+=\bar{z}^++d_o$ interconnected with the plant controller $v=\bar z^++g\!\left(\max\{|e_p|,|e_o|\}\right)\vartheta$. See Figure \ref{fig:feedback_optimization} for a block diagram interpretation. 

We show that the resulting feedback-optimization loop
satisfies the hypotheses of
Theorem~\ref{thm:gain_fxt_small_gain} and is FxT-ISS with
respect to $\mathcal A:=\{(y,z):y=z=y^\star\}$. In particular, define $q:=g\left(\max\{|e_p|,|e_o|\}\right)$, $u_p:=d_p-d_o$, and $u_o:=d_o$. The resulting error dynamics are
\begin{subequations}\label{eq:fo_interconnected_errors}
\begin{align}
e_p^+
&=q\vartheta+u_p,
\label{eq:fo_plant_error}\\
e_o^+
&=
G\left(G^{-1}(q\vartheta)+u_o\right).
\label{eq:fo_optimizer_error}
\end{align}
\end{subequations}
To verify Assumption~\ref{ass:sg_subsystems}, define $V_p(e_p):=|e_p|$ and $V_o(e_o):=|e_o|$. It is easy to see that the comparison bounds
\eqref{eq:sg_subsystem_bounds} hold with the identity function. Let $\psi_p(s)=\psi_o(s):=\psi(s):=s-g(s)$. Since $g$ is nonnegative and nondecreasing, $[s-\psi(s)]_+=g(s)$ and
\begin{align}
h_\psi(s)
&=
\max_{0\leq r\leq s}[r-\psi(r)]_+
=
\max_{0\leq r\leq s}g(r)
=g(s).
\label{eq:fo_hpsi}
\end{align}
Since $g^5(s)=0$ for all $s\geq0$, it follows that $\psi_p,\psi_o\in\mathcal E_{\mathrm{FxT}}^5$.

For the plant subsystem, using $|\vartheta|=1$ and the
monotonicity of $g$, we obtain
\begin{align}
V_p(e_p^+)
&\leq q+|u_p|
\nonumber\\
&=
g\left(
\max\{V_p(e_p),V_o(e_o)\}
\right)+|u_p|
\nonumber\\
&=
\max\left\{
g(V_p(e_p)),g(V_o(e_o))
\right\}+|u_p|.
\label{eq:fo_plant_max_bound}
\end{align}
Subtracting $V_p(e_p)$ and using
$g(s)=[s-\psi_p(s)]_+$ yields
\begin{align}
&V_p(e_p^+)-V_p(e_p)
\notag\\
&\quad\leq
\max\Big\{
-\min\{\psi_p(V_p(e_p)),V_p(e_p)\},
\notag\\
&\hspace{32mm}
g(V_o(e_o))-V_p(e_p)
\Big\}
+|u_p|.
\label{eq:fo_plant_increment}
\end{align}

For the optimizer subsystem, the global
$L$-Lipschitz continuity of $G$ gives
\begin{align}
V_o(e_o^+)
&=
\left|
G\left(G^{-1}(q\vartheta)+u_o\right)
\right|\leq
\left|
G\left(G^{-1}(q\vartheta)\right)
\right|
+L|u_o|
\nonumber\\
&=
q+L|u_o|
\nonumber\\
&=
\max\left\{
g(V_o(e_o)),g(V_p(e_p))
\right\}
+L|u_o|.
\label{eq:fo_optimizer_max_bound}
\end{align}
Consequently,
\begin{align}
&V_o(e_o^+)-V_o(e_o)
\notag\\
&\quad\leq
\max\Big\{
-\min\{\psi_o(V_o(e_o)),V_o(e_o)\},
\notag\\
&\hspace{32mm}
g(V_p(e_p))-V_o(e_o)
\Big\}
+L|u_o|.
\label{eq:fo_optimizer_increment}
\end{align}
Therefore, Assumption~\ref{ass:sg_subsystems} holds with $\gamma_{po}(s)=g(s)$, $\gamma_{op}(s)=g(s)$, $\chi_p(s)=s$, $\chi_o(s)=Ls$, $\sigma_p(s)=s$, and $\sigma_o(s)=s$. We next verify Assumption~\ref{ass:sg_fxt_condition}. Since
$\sigma_p=\sigma_o=\operatorname{id}$, the scaled internal
gains satisfy $b_{po}(s)=\sigma_p^{-1}\circ\gamma_{po}\circ\sigma_o(s)=g(s)$ and $b_{op}(s)=
\sigma_o^{-1}\circ\gamma_{op}\circ\sigma_p(s)
=g(s)$. Both functions are continuous, nonnegative, nondecreasing,
globally bounded by $3$, and vanish on $[0,1]$. Moreover, $b_{po}(s)=b_{op}(s)=g(s)<s$, for all $s>0$. Thus, Assumption~\ref{ass:sg_fxt_condition} holds. By Theorem \ref{thm:gain_fxt_small_gain}, the origin is FxT-ISS for the error dynamics \eqref{eq:fo_interconnected_errors}.


Finally, strong convexity and Lipschitz continuity of
$G=\nabla f$ imply that  $\mu|z-y^\star|\leq|G(z)|=|e_o|\leq L|z-y^\star|$. Since $e_p=y-z$, there exist constants $c_1,c_2>0$ such that $c_1 |(y,z)|_{\mathcal A} \leq |(e_p,e_o)| \leq c_2 |(y,z)|_{\mathcal A}$. Moreover, $(u_p,u_o)=(d_p-d_o,d_o)$ is a bounded linear transformation of $(d_p,d_o)$. Therefore, the feedback-optimization loop is FxT-ISS relative to $\mathcal A$ and with respect to the disturbances $(d_p,d_o)$. \hfill $\square$

\section{Proofs}
\label{secproofs}
%
%
\subsection{Proof of Lemma \ref{lem:general_fxt_extinction}:} We first prove sufficiency. Suppose that
conditions (a) and (b) hold. By condition (a) we have $g(s)=0$, for all $0\leq s\leq\varepsilon$. Moreover, condition (b) implies $0\leq g(s)\leq H$ for all $s\geq0$. Consequently, the function $h_{\psi}$ is continuous, nondecreasing, and satisfies $h_{\psi}(s)=0$ for all $0\leq s\leq\varepsilon$, and $h_{\psi}(s)\leq H$ for all $s\geq 0$.  We next show that $h_{\psi}(s)<s$ for all $s>0$. Indeed, since $\psi$ is positive definite, $g(r)<r$ for every $r>0$. Also, for any fixed $s>0$, continuity of $g$ guarantees the existence of $r_s\in[0,s]$ such that $h_{\psi}(s)=g(r_s)$. If $r_s=0$, then $h_{\psi}(s)=0<s$. Otherwise, $ h_{\psi}(s)=g(r_s)<r_s\leq s$ which establishes the desired bound. Now, suppose that $H\leq\varepsilon$. Then, $h_{\psi}^2(s)=0$ for all $s\geq 0$. Thus, $\psi\in\mathcal E_{\mathrm{FxT}}^2$. On the other hand, suppose that $H>\varepsilon$. By continuity and the properties of $h_{\psi}$, we have $\delta:=\min_{\varepsilon\leq s\leq H}\bigl(s-h_{\psi}(s)\bigr)>0$. Consider the scalar recursion  $v_{k+1}=h_{\psi}(v_k)$, with $v_0\geq 0$. By the properties of $h_{\psi}$, we have that
$v_1\leq H$. Furthermore, whenever
$v_k\in[\varepsilon,H]$, the iteration satisfies $v_{k+1}\leq v_k-\delta$. Therefore, after at most $L:=\left\lfloor\frac{H-\varepsilon}{\delta}\right\rfloor+1$ additional iterations, the trajectory enters $[0,\varepsilon)$. One further iteration then yields zero. Hence,  $h_{\psi}^M(s)=0$ for $s\geq 0$, with $M=
    3+
    \left\lfloor
        \frac{H-\varepsilon}{\delta}
    \right\rfloor$, which establishes sufficiency.

To prove necessity, suppose that $\psi\in\mathcal E_{\mathrm{FxT}}^M$ for some finite $M\geq 1$, so that $ h_{\psi}^M(s)=0$, for all $s\geq0$.
First, note that there must exist $\varepsilon>0$
such that $h_{\psi}(s)=0$, for all $0\leq s\leq\varepsilon$. Indeed, if $h_{\psi}(s)>0$ for every $s>0$,
then every finite iterate of $h$ would
remain strictly positive for every
positive initial condition, contradicting the fact that $ h_{\psi}^M(s)=0$, for all $s\geq0$. Thus, $h_{\psi}(\varepsilon)=0$ for some
$\varepsilon>0$. Since $h$ is
nonnegative and nondecreasing, the inequality follows. Since $g(s)\leq h_{\psi}(s)$, we have $g(s)=0$ for all  $0\leq s\leq\varepsilon$. Equivalently, $ \psi(s)\geq s$ for all $0\leq s\leq\varepsilon$.

Finally, suppose, by contradiction, that $h_{\psi}$ is unbounded. Since $h_{\psi}$ is continuous, nondecreasing, and satisfies $h_{\psi}(0)=0$, its range is $\mathbb R_{\geq 0}$. It follows inductively that every finite iterate $h_{\psi}^j$ is unbounded. In particular, $h_{\psi}^M$ is unbounded, contradicting the extinction property. Consequently, there exists $H<\infty$
such that $0\leq g(s)\leq h_{\psi}(s)\leq H$ for all $s\geq 0$. By the definition of $g$, we have  $[s-\psi(s)]_+\leq H$, which implies $ \psi(s)\geq s-H$ for all $s\geq 0$. \hfill $\blacksquare$
\subsection{Proof of Proposition \ref{prop:implication_to_dissipation}}
Let $\delta(s):=\min\{\psi(s),s\}$ for all $ s\geq 0$. For each $r\geq 0$, introduce the sublevel set $\Omega_r:=
    \{x\in\mathbb R^n:V(x)\leq\phi(r)\}$ and the function $q(r):=\max_{\substack{x\in\Omega_r\\|u|\leq r}}V(g(x,u))$. By \eqref{eq:general_lyap_bounds}, we have that $x\in\Omega_r$ implies $|x|_{\mathcal A}\leq\alpha_1^{-1}(\phi(r))$. Since $\mathcal A$ is compact, $\Omega_r$ is compact. Continuity of $g$ and $V$
therefore guarantees that $q(r)$ is well-defined and finite. Moreover, setting $u=0$ in \eqref{eq:prop_implication} gives $ V(g(x,0))
    \leq V(x)-\delta(V(x))$. In particular, $g(\mathcal A,0)\subseteq\mathcal A$, and hence
$q(0)=0$. We next establish that $\lim_{r\to 0^+}q(r)=0$. Suppose, by contradiction, that there exist $\varepsilon>0$, a sequence $r_j\to 0^+$, and pairs $(x_j,u_j)$ satisfying $V(x_j)\leq\phi(r_j)$ and
    $\|u_j\|\leq r_j$ such that $ V(g(x_j,u_j))\geq\varepsilon$. By \eqref{eq:general_lyap_bounds}, we have that $|x_j|_{\mathcal A}
    \leq\alpha_1^{-1}(\phi(r_j))\longrightarrow 0$. Compactness of $\mathcal A$ allows us to extract a subsequence, not relabeled, such that $x_j\to x_*\in\mathcal A$. Since $u_j\to 0$, continuity of $g$ and $V$ yields $ V(g(x_j,u_j))\longrightarrow V(g(x_*,0))=0,$ which is a contradiction. Thus, $\lim_{r\to 0^+}q(r)=0$ holds. Consequently, there exists $\bar q\in\mathcal K_\infty$ such that $q(r)\leq\bar q(r)$ for all  $r\geq 0$.

Next, let $\bar\delta(s):=\max_{0\leq t\leq s}\delta(t)$ and $\rho(r):=\bar q(r)+\bar\delta(\phi(r))$. Since $\delta$ is continuous and
positive definite, $\bar\delta$ is
continuous, nondecreasing, and satisfies
$\bar\delta(0)=0$. It follows that
$\rho\in\mathcal K_\infty$. We now distinguish two cases.

\emph{Case 1:} Suppose that $V(x)\geq\phi(|u|)$. Then, by \eqref{eq:prop_implication}, we have that  $V(g(x,u))-V(x)\leq-\delta(V(x))\leq-\delta(V(x))+\rho(|u|)$.
\emph{Case 2:} Suppose that $V(x)<\phi(|u|)$. Then, using the definition of $q$, as well as $\bar{q}$, we get $V(g(x,u))\leq q(|u|)\leq\bar q(|u|)$. Furthermore, $\delta(V(x))\leq\bar\delta(\phi(|u|))$. Therefore,
\begin{align}
    V(g(x,u))-V(x)
    &\leq\bar q(|u|)-V(x)\leq-\delta(V(x))
          +\bar q(|u|)
    \nonumber\\
    &\leq-\delta(V(x))+\rho(|u|),
    \label{eq:case_below_threshold}
\end{align}
where the second inequality follows
from $\delta(s)\leq s$. Combining the above two cases, and using the definition of $q$ we obtain \eqref{eq:additive_dissipation}.
Since $\psi\in\mathcal E_{\mathrm{FxT}}^M$
is unchanged, this inequality, together
with \eqref{eq:general_lyap_bounds}, establishes that
$V$ is an Extinction-Based FxT-ISS Lyapunov function
according to Definition~3.
\subsection{Proof of Lemma \ref{prop:newton_implication_impossibility}}
Suppose, by contradiction, that such functions exist. Since $\psi\in\mathcal E_{\mathrm{FxT}}^M$,
the associated map $ h(s):=\max_{0\leq r\leq s}[r-\psi(r)]_+$ satisfies $h^M(s)=0$ for every $s\geq0$, and by Lemma \ref{lem:general_fxt_extinction} there exists
$\varepsilon>0$ such that $\psi(s)\geq s$ for all $0\leq s\leq\varepsilon$. 

Let $\bar x\neq x^\star$ be
sufficiently close to $x^\star$, such that $\alpha_2(|\bar x-x^\star|)
 <\varepsilon$. The comparison bounds then yield $0<\alpha_1(|\bar x-x^\star|)\leq V(\bar x)<\varepsilon$. Notice that this conclusion does not
require continuity of $V$. Since $\phi\in\mathcal K_\infty$,
we can choose a nonzero disturbance
$\bar u$ sufficiently small such that $0<\phi(|\bar u|)<V(\bar x)$. Thus, the premise of the implication in \eqref{eq:prop_implication} holds at $(\bar x,\bar u)$. Moreover, the property that  $\psi(s)\geq s$ for all $0\leq s\leq\varepsilon$
and the property that $V(\bar{x})<\varepsilon$ imply $ \min\{\psi(V(\bar x)),V(\bar x)\}=V(\bar x)$.  Applying \eqref{eq:prop_implication}, we obtain
    $V(\bar x^+)-V(\bar x)\leq -V(\bar x)$ and therefore $V(\bar x^+)=0$. However, the perturbed Newton iteration \eqref{eq:disturbed_newton}
gives $ \bar x^+=x^\star+\bar u.$ Consequently, the lower comparison bound $V$ implies $ V(\bar x^+)\geq\alpha_1(|\bar u|)>0$, which is a contradiction. Therefore, no function $V$
satisfying \eqref{eq:general_lyap_bounds}, whether continuous or discontinuous, can satisfy the implication condition \eqref{eq:prop_implication}. \hfill $\blacksquare$
\subsection{Proof of Proposition \ref{prop:implication_fxt_iss}}
Define  $\delta(s):=\min\{\psi(s),s\}$. We distinguish two cases. \emph{Case 1:} If $V(x)\geq\phi(|u|)$, then by \eqref{eq:implication_condition},
\begin{align*}
V(x^+)-V(x)
&\leq-\delta(V(x))\leq-\delta(V(x))+\kappa(|u|).
\end{align*}
\emph{Case 2:} If $V(x)<\phi(|u|)$, then
\eqref{eq:subthreshold_condition} gives
$V(x^+)\leq\kappa(|u|)$. Since $\delta(s)\leq s$ for every $s\geq0$,
we obtain
\begin{align*}
V(x^+)-V(x)
&\leq\kappa(|u|)-V(x)\leq-\delta(V(x))+\kappa(|u|).
\end{align*}
Combining the above two cases
yields inequality \eqref{eq:additive_dissipation} with $\rho=\kappa$, which establishes the result.  \hfill $\blacksquare$
\subsection{Proof of Theorem \ref{thm:general_fxtiss}}
Let $v_k:=V(x(k))$. By \eqref{eq:additive_dissipation}, the nonnegativity of
$V$, and the definition of $h$ in \eqref{eq:general_h}, we have
\begin{equation}\label{iteration90}
v_{k+1}\leq h(v_k)+\rho(|u(k)|).
\end{equation}
We divide the proof into three steps.

\emph{Step 1:} By construction, $h$ is continuous, nondecreasing,
satisfies $h(0)=0$, and is globally bounded by some
constant $H>0$. For each $r\geq 0$, define
$T_r(s):=h(s)+\rho(r)$ so that whenever $|u|\leq r$,
the iteration \eqref{iteration90} becomes
$v_{k+1}\leq T_r(v_k)$.
Since $T_r(\cdot)$ is nondecreasing, consecutive
applications of $T_r(\cdot)$ yield
$v_{k+M}\leq T_r^M(v_k)$, provided
$\|u(j)\|\leq r$ for $k\leq j<k+M$.
Let
$\Gamma(r):=\sup_{s\geq 0}T_r^M(s)$.
Because $h(s)\leq H$, we have
$T_r(s)\leq H+\rho(r)$.
Monotonicity of $T_r$ therefore gives $\Gamma(r)\leq T_r^{M-1}(H+\rho(r))<\infty$. We claim that $\lim_{r\to 0^+}\Gamma(r)=0$. To see this, let
$H_*:=\sup_{s\geq 0}h(s)<\infty$.
By monotonicity and continuity of $h^{M-1}$,
we have
\begin{equation}\label{eqrrqe}
h^{M-1}(H_*)
=h^{M-1}\left(\sup_{s\geq 0}h(s)\right)
=\sup_{s\geq 0}h^M(s)=0.
\end{equation}
Moreover, monotonicity of $T_r$ yields $0\leq\Gamma(r)\leq T_r^{M-1}(H_*+\rho(r))$. Since $\rho(r)\to 0$ as $r\to 0^+$ and
$T_r=h+\rho(r)$, continuity of $h$ and
finite compositions, together with \eqref{eqrrqe}, gives
\[
\lim_{r\to 0^+}
T_r^{M-1}(H_*+\rho(r))
=h^{M-1}(H_*)=0.
\]
Consequently, $\lim_{r\to 0^+}\Gamma(r)=0$.
Since $\Gamma$ is finite and nondecreasing,
there exists $\overline{\Gamma}\in\mathcal{K}_\infty$
such that
$\Gamma(r)\leq\overline{\Gamma}(r)$
for all $r\geq 0$.

\emph{Step 2:} For $0\leq j<M$, define $d_j(r):=\sup_{s\geq 0}\left[T_r^j(s)-h^j(s)\right]$, where $T_r^0=h^0=\mathrm{id}$ is the identity map,
so that $d_0(r)=0$.
Since $T_r\geq h$ and both maps are nondecreasing,
the difference in $d_j$ is nonnegative.
We claim that each $d_j$ is finite,
nondecreasing, and satisfies
$\lim_{r\to 0^+}d_j(r)=0$.

To see this, let $r_0>0$ and define
$B:=H+\rho(r_0)$. Then, for every
$r\in[0,r_0]$ and $s\geq 0$, we have
$0\leq T_r(s)\leq B$ and
$0\leq h(s)\leq B$.
Thus, all iterates of $T_r$ and $h$,
after the first, belong to the common compact
interval $[0,B]$. By continuity, $h$ is
uniformly continuous on $[0,B]$. Let $\omega_B(\delta):=
\sup_{\substack{x,y\in[0,B]\\|x-y|\leq\delta}}
|h(x)-h(y)|$, which satisfies
$\omega_B(\delta)\to 0$ as $\delta\to 0^+$.
For $j=1$, we have $d_1(r)=\rho(r)$.
Moreover, for every $j\geq 1$ and $s\geq 0$,
\begin{align*}
|T_r^{j+1}(s)-h^{j+1}(s)|&=|h(T_r^j(s))+\rho(r)-h(h^j(s))|\\
&\leq\rho(r)+\omega_B\bigl(|T_r^j(s)-h^j(s)|\bigr).
\end{align*}
Taking the supremum over $s\geq 0$ and using
the monotonicity of $\omega_B$ gives $d_{j+1}(r)\leq\rho(r)+\omega_B(d_j(r))$. Since $\rho(r)\to 0$ as $r\to 0^+$,
we have $d_1(r)\to 0$.
Suppose that $d_j(r)\to 0$ for some fixed
$j\geq 1$. The preceding inequality and
the continuity of $\omega_B$ at zero imply
$\lim_{r\to 0^+}d_{j+1}(r)=0$.
The result follows by induction.
Finiteness and monotonicity follow from
the preceding bounds and the monotonicity
of $T_r$ with respect to $r$.

\emph{Step 3:}
Since $M$ is finite, the function $d(r):=\max_{0\leq j<M}d_j(r)$ is finite, nondecreasing, and satisfies
$d(r)\to 0$ as $r\to 0^+$.
Choose $\overline d\in\mathcal{K}_\infty$
such that $d(r)\leq\overline d(r)$
for all $r\geq 0$.
Let
$U_k:=\max_{0\leq j<k}|u(j)|$,
with $U_0:=0$.
Then, for $0\leq k<M$, repeated application
of \eqref{iteration90} and monotonicity of $T_{U_k}(\cdot)$
yield
\begin{equation}\label{boundvk}
v_k
\leq T_{U_k}^k(v_0)
\leq h^k(v_0)+\overline d(U_k)
\leq\max\{2h^k(v_0),2\overline d(U_k)\}.
\end{equation}
For $k\geq M$, Step 1 gives $v_k\leq\Gamma\left(\max_{k-M\leq j<k}|u(j)|\right)\leq\overline\Gamma(U_k)$. Define
$\gamma(r):=\alpha_1^{-1}(\widetilde\gamma(r))$,
where $\widetilde\gamma(r):=\max\{2\overline d(r),\overline\Gamma(r)\}$, and

\vspace{-0.2cm}
\begin{small}
\begin{equation}\label{klbound1}
\beta(s,k):=
\begin{cases}
\displaystyle
\alpha_1^{-1}\left(
2\Big[
\max_{k\leq j<M}h^j(\alpha_2(s))
+\alpha_2(s)
\Big]\right),
&0\leq k<M,\\
0,&k\geq M.
\end{cases}
\end{equation}
\end{small}
For every $k<M$, 
$\beta(\cdot,k)$ is continuous, strictly
increasing, and vanishes at zero.
Moreover, the maximization in \eqref{klbound1}
is taken over a decreasing family of
index sets as $k$ increases. Hence,
$\beta(s,k)$ is nonincreasing in $k$,
and $\beta(s,k)=0,\qquad k\geq M$. Thus, $\beta\in\mathcal{KL}_{\mathrm{FxT}}$.
Combining \eqref{boundvk}, the bound for $k\geq M$,
and the comparison bounds \eqref{eq:general_lyap_bounds}, we obtain  $|x(k)|_{\mathcal A}
\leq
\max\left\{
\beta(|x(0)|_{\mathcal A},k),
\gamma(U_k)
\right\}$, $\forall~k\geq 0$. \hfill $\blacksquare$
\subsection{Proof of Theorem \ref{thm:gain_fxt_small_gain}}
For each $i\in\{1,2\}$, let $Y_i(x_i):=\sigma_i^{-1}(V_i(x_i))$,
%
and consider the composite Lyapunov candidate $W(x):=\max_{i\in\{1,2\}}Y_i(x_i)$.
We divide the proof into three main steps.

\emph{Step 1:}
Let $\underline{\alpha}(s)=
\min_{i\in\{1,2\}}\{\sigma_i^{-1}\circ\underline{\alpha}_i\left(\frac{s}{\sqrt{2}}\right)
\}$ and $\overline{\alpha}(s):=
\max_{i\in\{1,2\}}
\left\{
\sigma_i^{-1}\circ\overline{\alpha}_i(s)
\right\}$. Both $\underline{\alpha}$ and $\overline{\alpha}$ belong to
$\mathcal K_\infty$, and satisfy $\underline{\alpha}(|x|)\leq W(x)\leq\overline{\alpha}(|x|)$.
%
Since $h_{\psi_i}(s)=
\max_{0\leq r\leq s}[r-\psi_i(r)]_+$, we have that $[s-\psi_i(s)]_+
\leq h_{\psi_i}(s)$. Therefore, \eqref{eq:sg_subsystem_max_form} implies
\begin{align*}
V_i(x_i^+)
\leq
\max\Big\{
&h_{\psi_i}(V_i(x_i)),
\gamma_{ij}(V_j(x_j))
\Big\}
+\chi_i(|u_i|).
\end{align*}

\emph{Step 2:} Let  $H(s):=\max\{a_1(s),a_2(s),b_{12}(s),b_{21}(s)\}$. Because
$\psi_i\in\mathcal E_{\mathrm{FxT}}^{M_i}$,
the function $h_{\psi_i}$ is continuous, nonnegative,
nondecreasing, globally bounded, vanishes on a neighborhood
of the origin, and satisfies $h_{\psi_i}(s)<s$ for all $s>0$. These properties are preserved under conjugation by
$\sigma_i\in\mathcal K_\infty$. Thus, each $a_i$ is
continuous, nonnegative, nondecreasing, globally bounded,
vanishes on a neighborhood of the origin, and satisfies $a_i(s)<s$ for all $s>0$. By the definition of $b_{ij}$ and the properties of
$\gamma_{ij}$ and $\sigma_i$, each function $b_{ij}$ is
continuous, nonnegative, and nondecreasing, with
$b_{ij}(0)=0$. Assumption~\ref{ass:sg_fxt_condition}
guarantees that $b_{ij}$ is globally bounded and vanishes
on a neighborhood of the origin, whereas
\eqref{eq:sg_internal_gain_condition} gives $b_{ij}(s)<s$, $\forall s>0$. Consequently, $H$ is continuous,
nonnegative, nondecreasing, globally bounded, vanishes on a
neighborhood of the origin, and satisfies $H(s)<s$ for all $s>0$. Let $\psi(s):=s-H(s)$, which is continuous and positive definite. Moreover, because $H$ is nondecreasing,
\begin{align}
h_\psi(s)
&=
\max_{0\leq r\leq s}[r-\psi(r)]_+=
\max_{0\leq r\leq s}H(r)
=H(s).
\label{eq:sg_hpsi_H}
\end{align}
The boundedness, dead-zone, and strict-contraction properties
of $H$, together with
Lemma~\ref{lem:general_fxt_extinction}, imply that there
exists $M\in\mathbb Z_{\geq1}$ such that $H^M(s)=0$, i.e., $\psi\in\mathcal E_{\mathrm{FxT}}^M$.

\emph{Step 3:} Since the functions $a_i$ and $b_{ij}$ are globally bounded, there exists $B>0$ such that $a_i(s)\leq B$, and $b_{ij}(s)\leq B$ for every $s\geq0$, $i\in\{1,2\}$, and $j\neq i$. For each $i\in\{1,2\}$, define
\begin{equation}\label{eq:sg_transformed_input_gain}
\widehat{\chi}_i(r)
:=
\max_{0\leq t\leq\sigma_i(B)}
\left[
\sigma_i^{-1}(t+\chi_i(r))
-\sigma_i^{-1}(t)
\right].
\end{equation}
The function $\widehat{\chi}_i$ belongs to
$\mathcal K_\infty$. For every $i\in\{1,2\}$ and $j\neq i$, let $A_i(x):=
\max\left\{h_{\psi_i}(V_i(x_i)),\gamma_{ij}(V_j(x_j))\right\}$. Since $V_i(x_i)=\sigma_i(Y_i(x_i))$ and
$Y_i(x_i)\leq W(x)$, monotonicity gives
\begin{align*}
\sigma_i^{-1}(A_i(x))
&\leq
\max\left\{
a_i(W(x)),b_{ij}(W(x))
\right\}\leq H(W(x)).
\end{align*}
Thus, $\sigma_i^{-1}(A_i(x))\leq B$, and $A_i(x)\leq\sigma_i(B)$. Applying $\sigma_i^{-1}$ to
\eqref{eq:sg_enveloped_subsystem} and using
\eqref{eq:sg_transformed_input_gain} yields
\begin{align*}
Y_i(x_i^+)&=\sigma_i^{-1}(V_i(x_i^+))\leq
\sigma_i^{-1}\left(
A_i(x)+\chi_i(|u_i|)
\right)
\nonumber\\
&\leq
\sigma_i^{-1}(A_i(x))
+\widehat{\chi}_i(|u_i|)\leq
H(W(x))
+\widehat{\chi}_i(|u_i|).
\end{align*}
Let $\rho(r):=\max_{i\in\{1,2\}}\widehat{\chi}_i(r)$, and note that $\rho\in\mathcal K_\infty$. Since
$|u_i|\leq|u|$, taking the maximum over
$i\in\{1,2\}$ in $Y_i(x_i^+)$ gives
\begin{equation}\label{eq:sg_composite_update}
W(x^+)
\leq
H(W(x))+\rho(|u|).
\end{equation}
Subtracting $W(x)$ and using
$\psi(s)=s-H(s)$ gives
\begin{align}
W(x^+)-W(x)
&\leq
-\psi(W(x))+\rho(|u|)
\nonumber\\
&=
-\min\{\psi(W(x)),W(x)\}
+\rho(|u|),
\label{eq:sg_composite_dissipation}
\end{align}
where the equality follows from
$0\leq\psi(s)\leq s$. Therefore, $W$ is a Extinction-Based
FxT-ISS Lyapunov function for the interconnected system, and
by Theorem~\ref{thm:general_fxtiss}, it is FxT-ISS. \hfill $\blacksquare$
\section{CONCLUSIONS}
\label{section:conclusions}
This paper developed a framework for the study of fixed-time input-to-state stability in discrete-time systems. Using a class of dissipation functions characterized by a fixed-time extinction property, we established Lyapunov-based sufficient conditions for FxT-ISS that extend beyond commonly used two-power constructions. Building on these results, we developed a small-gain theorem for interconnected FxT-ISS systems, providing conditions under which fixed-time stability is preserved in the presence of interconnections. 

\bibliographystyle{IEEEtran}
\bibliography{Biblio.bib}

\end{document}